\input amstex
\documentstyle{amsppt}
%
\catcode`@=11
\redefine\output@{%
  \def\break{\penalty-\@M}\let\par\endgraf
  \ifodd\pageno\global\hoffset=105pt\else\global\hoffset=8pt\fi  
  \shipout\vbox{%
    \ifplain@
      \let\makeheadline\relax \let\makefootline\relax
    \else
      \iffirstpage@ \global\firstpage@false
        \let\rightheadline\frheadline
        \let\leftheadline\flheadline
      \else
        \ifrunheads@ 
        \else \let\makeheadline\relax
        \fi
      \fi
    \fi
    \makeheadline \pagebody \makefootline}%
  \advancepageno \ifnum\outputpenalty>-\@MM\else\dosupereject\fi
}
\def\Beta{\mathchar"0\hexnumber@\rmfam 42}
\catcode`\@=\active
\nopagenumbers
\chardef\textvolna='176

\chardef\bigalpha='013
\def\negskp{\hskip -2pt}

\def\Img{\operatorname{Im}}

\chardef\degree="5E
\def\const{\operatorname{const}}
\def\blue#1{#1}

\gdef\darkred#1{#1}
\catcode`#=11\def\diez{#}\catcode`#=6
\catcode`&=11\catcode`&=4
\catcode`_=11\def\podcherkivanie{_}\catcode`_=8
\catcode`\^=11\catcode`\^=7
\catcode`~=11\def\volna{~}\catcode`~=\active
\def\mycite#1{\cite{\blue{#1}}\immediate\special{ps:
     ShrHPSdict begin /ShrBORDERthickness 0 def}}
\def\myciterange#1#2#3#4{\cite{\blue{#2#3#4}}\immediate\special{ps:
     ShrHPSdict begin /ShrBORDERthickness 0 def}}
\def\mytag#1{%
    \tag#1}
\def\mythetag#1{\thetag{\blue{#1}}\immediate\special{ps:
     ShrHPSdict begin /ShrBORDERthickness 0 def}}
\def\myrefno#1{\no#1}
\def\myhref#1#2{\blue{#2}\immediate\special{ps:
     ShrHPSdict begin /ShrBORDERthickness 0 def}}
\def\myEarXivlink{\myhref{http://arXiv.org}{http:/\negskp/arXiv.org}}

\def\mytheorem#1{\csname proclaim\endcsname{Theorem #1}}
\def\mytheoremwithtitle#1#2{\csname proclaim\endcsname{Theorem #1#2}}
\def\mythetheorem#1{\blue{#1}\immediate\special{ps:
     ShrHPSdict begin /ShrBORDERthickness 0 def}}
\def\mylemma#1{\csname proclaim\endcsname{Lemma #1}}
\def\mylemmawithtitle#1#2{\csname proclaim\endcsname{Lemma #1#2}}
\def\mythelemma#1{\blue{#1}\immediate\special{ps:
     ShrHPSdict begin /ShrBORDERthickness 0 def}}
\def\mycorollary#1{\csname proclaim\endcsname{Corollary #1}}

\def\myconjecture#1{\csname proclaim\endcsname{Conjecture #1}}
\def\myconjecturewithtitle#1#2{\csname proclaim\endcsname{Conjecture #1#2}}
\def\mytheconjecture#1{\blue{#1}\immediate\special{ps:
     ShrHPSdict begin /ShrBORDERthickness 0 def}}
\def\myproblem#1{\csname proclaim\endcsname{Problem #1}}
\def\myproblemwithtitle#1#2{\csname proclaim\endcsname{Problem #1#2}}


\pagewidth{360pt}
\pageheight{606pt}
\topmatter
\title
Asymptotic estimates for roots of the cuboid characteristic 
equation in the linear region.
\endtitle
\rightheadtext{Asymptotic estimates for roots \dots}
\author
Ruslan Sharipov
\endauthor
\address Bashkir State University, 32 Zaki Validi street, 450074 Ufa, Russia
\endaddress
\email
\myhref{mailto:r-sharipov\@mail.ru}{r-sharipov\@mail.ru}
\endemail
\abstract
     A perfect cuboid is a rectangular parallelepiped whose edges, whose face 
diagonals, and whose space diagonal are of integer lengths. The second cuboid 
conjecture specifies a subclass of perfect cuboids described by one Diophantine 
equation of tenth degree and claims their non-existence within this subclass. 
This Diophantine equation is called the cuboid characteristic equation. It has 
two parameters. The linear region is a domain on the coordinate plane of these 
two parameters given by certain linear inequalities. In the present paper 
asymptotic expansions and estimates for roots of the characteristic equation 
are obtained in the case where both parameters tend to infinity staying within 
the linear region. Their applications to the cuboid problem are discussed. 
\endabstract
\endtopmatter
\TagsOnRight
\document

\head
1. Introduction.
\endhead
     The cuboid characteristic equation in the case of the second cuboid conjecture 
is a polynomial Diophantine equation with two parameters $p$ and $q$:
$$
\hskip -2em
Q_{p\kern 0.6pt q}(t)=0.
\mytag{1.1}
$$
The polynomial $Q_{p\kern 0.6pt q}(t)$ from \mythetag{1.1} is given
by the explicit formula 
$$
\gathered
Q_{p\kern 0.6pt q}(t)=t^{10}+(2\,q^{\kern 0.7pt 2}+p^{\kern 1pt 2})\,(3
\,q^{\kern 0.7pt 2}-2\,p^{\kern 1pt 2})\,t^8+(q^{\kern 0.5pt 8}+10
\,p^{\kern 1pt 2}\,q^{\kern 0.5pt 6}+\\
+\,4\,p^{\kern 1pt 4}\,q^{\kern 0.5pt 4}-14\,p^{\kern 1pt 6}\,q^{\kern 0.7pt 2}
+p^{\kern 1pt 8})\,t^6-p^{\kern 1pt 2}\,q^{\kern 0.7pt 2}
\,(q^{\kern 0.5pt 8}-14\,p^{\kern 1pt 2}\,q^{\kern 0.5pt 6}
+4\,p^{\kern 1pt 4}\,q^{\kern 0.5pt 4}+\\
+\,10\,p^{\kern 1pt 6}\,q^{\kern 0.7pt 2}+p^{\kern 1pt 8})\,t^4
-p^{\kern 1pt 6}\,q^{\kern 0.5pt 6}\,(q^{\kern 0.7pt 2}
+2\,p^{\kern 1pt 2})\,(3\,p^{\kern 1pt 2}-2\,q^{\kern 0.7pt 2})\,t^2
-q^{\kern 0.7pt 10}\,p^{\kern 1pt 10}.
\endgathered\quad
\mytag{1.2}
$$
The tenth degree polynomial \mythetag{1.2} is related to the perfect cuboid 
problem through the following theorem (see Theorem~8.1 in \mycite{1} or in
\mycite{2}).
\mytheorem{1.1} A triple of positive integer numbers $p$, $q$, and $t$ satisfying 
the equation \mythetag{1.1} and such that $p\neq q$ are coprime produces a perfect
cuboid if and only if the following inequalities are
fulfilled: 
$$
\xalignat 4
&\hskip -2em
t>p^{\kern 1pt 2},
&&t>p\,q,
&&t>q^{\kern 0.7pt 2},
&&(p^{\kern 1pt 2}+t)\,(p\,q+t)>2\,t^2.
\qquad
\mytag{1.3}
\endxalignat
$$
\endproclaim
     Once a triple of numbers $p,\,q,\,t$ obeying Theorem~\mythetheorem{1.1} is 
found, there is a definite procedure for producing a perfect cuboid from them. 
\pagebreak First of all three rational numbers $\alpha$, $\beta$, and $\upsilon$ 
are produced in one of the two ways: 1) using the formulas
$$
\xalignat 3
&\hskip -2em
\alpha=\frac{p^{\kern 1pt 2}}{t},
&&\beta=\frac{p\,q}{t},
&&\upsilon=\frac{q^{\kern 0.7pt 2}}{t},
\mytag{1.4}
\endxalignat
$$
or 2) using the other three formulas 
$$
\xalignat 3
&\hskip -2em
\alpha=\frac{p\,q}{t},
&&\beta=\frac{p^{\kern 1pt 2}}{t},
&&\upsilon=\frac{q^{\kern 0.7pt 2}}{t}.
\mytag{1.5}
\endxalignat
$$
Both ways \mythetag{1.4} and \mythetag{1.5} are acceptable.\par
     Then the rational number $z$ is produced using the following
formula for it:
$$
\hskip -2em
z=\frac{(1+\upsilon^2)\,(1-\beta^2)\,(1+\alpha^2)}{2\,(1+\beta^2)\,(1
-\alpha^2\,\upsilon^2)}.
\mytag{1.6}
$$
And finally, the numbers $\alpha$, $\beta$, $\upsilon$ along with
$z$ from \mythetag{1.6} are used in the formulas 
$$
\xalignat 2
&\hskip -2em 
\frac{x_1}{L}=\frac{2\,\upsilon}{1+\upsilon^2},
&&\frac{d_1}{L}=\frac{1-\upsilon^2}{1+\upsilon^2},\\
\vspace{1ex}
&\hskip -2em
\frac{x_2}{L}=\frac{2\,z\,(1-\upsilon^2)}{(1+\upsilon^2)\,(1+z^2)},
&&\frac{x_3}{L}=\frac{(1-\upsilon^2)\,(1-z^2)}{(1+\upsilon^2)\,(1+z^2)},
\mytag{1.7}\\
\vspace{1ex}
&\hskip -2em
\frac{d_2}{L}=\frac{(1+\upsilon^2)\,(1+z^2)+2\,z(1-\upsilon^2)}
{(1+\upsilon^2)\,(1+z^2)}\,\beta,
&&\frac{d_3}{L}=\frac{2\,(\upsilon^2\,z^2+1)}{(1+\upsilon^2)\,(1+z^2)}\,\alpha.
\quad
\endxalignat
$$
They produce six rational numbers in the right hand sides of the formulas
\mythetag{1.7}. Then $L$ is chosen as a common denominator for all these 
six rational numbers. Such a choice assures that $x_1$, $x_2$, $x_3$, $d_1$, 
$d_2$, $d_3$ are integer numbers. They are edges and face diagonals of
a perfect cuboid, while $L$ is its space diagonal. In the other words,
the integer numbers $x_1$, $x_2$, $x_3$, $d_1$, $d_2$, $d_3$, and $L$
satisfy the cuboid equations 
$$
\xalignat 2
&\hskip -2em
x_1^2+x_2^2+x_3^2=L^2,
&&x_2^2+x_3^2=d_1^{\kern 1pt 2},\\
\vspace{-1.7ex}
\mytag{1.8}\\
\vspace{-1.7ex}
&\hskip -2em
x_3^2+x_1^2=d_2^{\kern 1pt 2},
&&x_1^2+x_2^2=d_3^{\kern 1pt 2}.
\endxalignat
$$\par
     The formulas \mythetag{1.7} were derived from \mythetag{1.8} in \mycite{3}
for the general case where, instead of \mythetag{1.2}, a twelfth degree polynomial
arises. It reduces to the polynomial \mythetag{1.2} in a special case, which is 
called the case of the second cuboid conjecture. The second cuboid conjecture
itself is formulated as follows (see \mycite{4}). 
\myconjecture{1.1} For any positive coprime integers $p\neq q$ the 
polynomial $Q_{p\kern 0.6pt q}(t)$ is irreducible over the ring of 
integer numbers. 
\endproclaim
    Conjecture~\mytheconjecture{1.1} means that the equation \mythetag{1.1}
has no integer roots. However, this conjecture is not yet proved nor disproved. 
It is just a conjecture. Therefore in the present paper we consider the roots
of the equation \mythetag{1.1} and study their dependence on $p$ and $q$. This
research continues the research from \myciterange{3}{3}{--}{7} and 
\myciterange{1}{1}{,\,}{2}. As for the Diophantine equations \mythetag{1.8},
they are being studied for almost 300 years. For the history and various 
approaches to them the reader is referred to \myciterange{8}{8}{--}{50}. The
approach of the papers \myciterange{51}{51}{--}{63} is based on so-called 
multisymmetric polynomials. It is different from the approach of the present
paper. Therefore we do not consider the papers \myciterange{51}{51}{--}{63}
below.\par
     The linear region associated with the cuboid characteristic polynomial 
\mythetag{1.2} in the case of the second cuboid conjecture was defined in
\mycite{2}. It is a domain in the positive quadrant of the 
$p\,q$\,-\,coordinate plane given by the linear inequalities
$$
\xalignat 2
&\frac{q}{59}<p,
&&p<59\,q.
\mytag{1.9}
\endxalignat
$$
The main goal of the present paper is to obtain asymptotic expansions and 
estimates for the roots of the characteristic equation \mythetag{1.1}
as $p\to+\infty$ and $q\to+\infty$ simultaneously staying within the linear
region \mythetag{1.9}.\par
\head
2. Asymptotic expansions with constant ratio. 
\endhead
    It is easy to see that the inequalities defining the linear 
region \mythetag{1.9} set upper and lower bounds for the ratio of two parameters 
$p$ and $q$. They are written as 
$$
\hskip -2em
\frac{1}{59}<\frac{p}{q}<59.
\mytag{2.1}
$$
Due to \mythetag{2.1} below we consider the case where
$$
\xalignat 3
&\hskip -2em
p\to+\infty,
&&q\to+\infty,
&&\frac{p}{q}\to\theta\neq\infty
\mytag{2.2}
\endxalignat
$$
and where $\theta$ is a rational number. In this case we can write
$$
\hskip -2em
\theta=\frac{a_{12}}{a_{22}},
\mytag{2.3}
$$
where $a_{11}$ and $a_{12}$ are two positive coprime integers, i\.\,e\.
the fraction \mythetag{2.3} is irreducible. For any two positive 
coprime integers $a_{11}$ and $a_{12}$ there are two other coprime 
integers $a_{21}$ and $a_{22}$ such that
$$
\hskip -2em
a_{11}\,a_{22}-a_{21}\,a_{12}=1. 
\mytag{2.4}
$$
This fact follows from the Euclidean division algorithm (see 
\mycite{64} or \mycite{65}).\par
     Using the numbers $a_{11}$, $a_{12}$, $a_{21}$, $a_{22}$ from
\mythetag{2.3} and \mythetag{2.4}, we define two matrices
$$
\xalignat 2
&\hskip -2em
S=\Vmatrix 
a_{11} & a_{12}\\
\vspace{1ex}
a_{21} & a_{22}
\endVmatrix,
&&T=\Vmatrix 
a_{22} & -a_{12}\\
\vspace{1ex}
-a_{21} & a_{11}
\endVmatrix.
\mytag{2.5}
\endxalignat 
$$
The equality \mythetag{2.4} means that $\det S=1$ and $\det T=1$. Moreover,
the matrices \mythetag{2.5} are inverse to each other. We use them as 
transition matrices (see \mycite{66}) and define the following change of 
coordinates in the $p\,q$\,-\,coordinate plane:
$$
\xalignat 2
&\hskip -2em
\Vmatrix 
\tilde p\\
\vspace{1ex}
\tilde q
\endVmatrix
=\Vmatrix 
a_{22} & -a_{12}\\
\vspace{1ex}
-a_{21} & a_{11}
\endVmatrix
\cdot
\Vmatrix 
p\\
\vspace{1ex}
q
\endVmatrix,
&&\Vmatrix 
p\\
\vspace{1ex}
q
\endVmatrix
=\Vmatrix 
a_{11} & a_{12}\\
\vspace{1ex}
a_{21} & a_{22}
\endVmatrix
\cdot
\Vmatrix 
\tilde p\\
\vspace{1ex}
\tilde q
\endVmatrix.
\qquad
\mytag{2.6}
\endxalignat 
$$
The formulas \mythetag{2.6} can be written in a non-matrix way:
$$
\xalignat 2
&\hskip -2em
\cases 
\tilde p=\hphantom{-}a_{22}\,p-a_{12}\,q,\\
\tilde q=-a_{21}\,p+a_{11}\,q,
\endcases
&&\cases 
p=a_{11}\,\tilde p+a_{12}\,\tilde q,\\
q=a_{21}\,\tilde p+a_{22}\,\tilde q.
\endcases
\quad
\mytag{2.7}
\endxalignat 
$$
Using \mythetag{2.7}, \mythetag{2.2}, \mythetag{2.3}, and \mythetag{2.4}, 
we derive 
$$
\align
&\hskip -2em
\frac{\tilde p}{q}=a_{22}\,\frac{p}{q}-a_{12}\to a_{22}\,\theta-a_{12}
=\frac{a_{22}\,a_{12}-a_{12}\,a_{22}}{a_{22}}=0,
\mytag{2.8}\\
&\hskip -2em
\frac{\tilde q}{q}=-a_{21}\,\frac{p}{q}+a_{11}\to -a_{21}\,\theta+a_{11}
=\frac{-a_{21}\,a_{12}+a_{11}\,a_{22}}{a_{22}}=\frac{1}{a_{22}}
\mytag{2.9}
\endalign
$$
as $q\to+\infty$. Relying on \mythetag{2.8} and \mythetag{2.9}, we use
the more restrictive condition 
$$
\hskip -2em
\tilde p=\const\text{\ \ as \ }\tilde q\to +\infty 
\mytag{2.10}
$$
when passing to the new variables $\tilde p$ and $\tilde q$. It is worth
to note the following lemma describing $\tilde p$ and $\tilde q$
in \mythetag{2.7}. 
\mylemma{2.1} If $p$ and $q$ are coprime then the numbers $\tilde p$ and 
$\tilde q$ produced according to the first couple of the formulas \mythetag{2.7} 
are also coprime.
\endproclaim
     The proof is immediate from the second couple of the formulas 
\mythetag{2.7}. Indeed, if $\tilde p$ and $\tilde q$ have some common divisor
$r$, then from the second couple of the formulas \mythetag{2.7} we derive 
that $r$ is a common divisor $p$ and $q$, which contradicts their coprimality.
\par
     Using \mythetag{2.7}, let's substitute $p=a_{11}\,\tilde p+a_{12}\,
\tilde q$ and $q=a_{21}\,\tilde p+a_{22}\,\tilde q$ into the polynomial 
\mythetag{1.2}. As a result we get another polynomial $Q_{\tilde p\kern 0.6pt 
\tilde q}(t)$. This polynomial is given by an explicit formula. However, the 
formula for the polynomial $Q_{\tilde p\kern 0.6pt \tilde q}(t)$ is rather huge.
It is placed to the ancillary file \darkred{{\tt strategy\kern -0.5pt\_\kern 1.5pt 
formulas\_\kern 0.5pt 03.txt}} in a machine-readable form.\par
     Using the polynomial $Q_{\tilde p\kern 0.6pt \tilde q}(t)$, we replace
the equation \mythetag{1.1} by the equation
$$
\hskip -2em
Q_{\tilde p\kern 0.6pt \tilde q}(t)=0.
\mytag{2.11}
$$
It is a tenth degree equation with respect to the variable $t$. Like \
$Q_{p\kern 0.6pt q}(t)$, the polynomial $Q_{\tilde p\kern 0.6pt \tilde q}(t)$
in \mythetag{2.11} is even with respect to $t$. Along with each root $t$ 
it has the opposite root $-t$. Therefore we use the condition
$$
\hskip -2em
\cases \text{$t>0$ \ if \ $t$ \ is a real root,}\\
\text{$\Img(t)>0$ \ if \ $t$ \ is a complex root,}
\endcases
\mytag{2.12}
$$
in order to divide the roots of the equation \mythetag{2.11} into two groups. 
We denote through $t_1$, $t_2$, $t_3$, $t_4$, $t_5$ those roots that obey the 
conditions \mythetag{2.12}. Then $t_6$, $t_7$, $t_8$, $t_9$, $t_{10}$ are 
opposite roots of the equation \mythetag{2.11}:
$$
\xalignat 5
&\hskip -1em
t_6=-t_1,&&t_7=-t_2,&&t_8=-t_3,&&t_9=-t_4,&&t_{10}=-t_5.
\qquad\quad
\mytag{2.13}
\endxalignat
$$\par
    Typically, asymptotic expansions for roots of a polynomial equation look 
like power series (see \mycite{67}). By analogy to (2.3) in \mycite{1} and
according to \mythetag{2.10} we write
$$
\hskip -2em
t_i(\tilde p,\tilde q)=C_i\,\tilde q^{\,\alpha_i}
\biggl(1+\sum^\infty_{s=1}\beta_{is}\,\tilde q^{-s}\biggr)
\text{\ \ as \ }\tilde q\to+\infty.
\mytag{2.14}
$$
The coefficients $C_i$ in \mythetag{2.14} should be nonzero: 
$C_i\neq 0$.\par
     The exponents $\alpha_i$ and the coefficients $C_i$ in \mythetag{2.14}
are determined using the Newton polygon (see \mycite{1}). The Newton polygon
associated with the polynomial $Q_{\tilde p\kern 0.6pt \tilde q}(t)$ in
\mythetag{2.11} is a triangle \vadjust{\vskip 536pt\hbox to 0pt{\kern 40pt
\includegraphics{Strategy03.eps}\hss}\vskip 0pt}(see Fig\.~2.1 below).  
Its boundary consists of three parts --- the upper part, the lower part, and
the vertical part. The upper part is drawn in green, the lower part is drawn 
in red. The upper part of the Newton polygon in Fig\.~2.1 is a segment of a 
straight line. It comprises six nodes associated with the equation 
\mythetag{2.11}. Here are the coefficients of these nodes: 
$$
\hskip -2em
\aligned
&A_{\kern 0.5pt 10\kern 2pt 0}=1,\\
&A_{\kern 0.5pt 8\kern 1.5pt 4}=6\,a_{22}^4-2\,a_{12}^4-a_{12}^2\,a_{22}^2,\\
&A_{\kern 0.5pt 6\kern 1.5pt 8}=10\,a_{12}^2\,a_{22}^6+a_{12}^8+a_{22}^8
+4\,a_{12}^4\,a_{22}^4-14\,a_{12}^6\,a_{22}^2,\\
&A_{\kern 0.5pt 4\kern 1.5pt 12}=14\,a_{12}^4\,a_{22}^8-4\,a_{12}^6\,a_{22}^6
-a_{12}^2\,a_{22}^{10}-10\,a_{12}^8\,a_{22}^4-a_{12}^{10}\,a_{22}^2,\\
&A_{\kern 0.5pt 2\kern 1.5pt 16}=a_{12}^8\,a_{22}^8-6\,a_{12}^{10}\,a_{22}^6
+2\,a_{12}^6\,a_{22}^{10},\\
&A_{\kern 0.5pt 0\kern 1.5pt 20}=-a_{12}^{10}\,a_{22}^{10}.
\endaligned
\mytag{2.15}
$$\par
     According to \mycite{1}, the exponent $\alpha_i$ in \mythetag{2.11} is 
determined by the slope of the upper boundary of the Newton triangle in 
Fig\.~2.1 by means of the formula $\alpha_i=-k$. In our case $k=-2$, hence for
$\alpha_i$ we derive 
$$
\hskip -2em
\alpha_i=2.
\mytag{2.16}
$$
The exponent \mythetag{2.16} is common for all roots of the equation 
\mythetag{2.11}. The coefficient $C_i$ in \mythetag{2.14} is determined
by the following equation:
$$
\hskip -2em
\aligned
A_{\kern 0.5pt 10\kern 2pt 0}\ {C_i}^{10}
&+A_{\kern 0.5pt 8\kern 1.5pt 4}\ {C_i}^8
+A_{\kern 0.5pt 6\kern 1.5pt 8}\ {C_i}^6\,+\\
&+\,A_{\kern 0.5pt 4\kern 1.5pt 12}\ {C_i}^4
+A_{\kern 0.5pt 2\kern 1.5pt 16}\ {C_i}^2
+A_{\kern 0.5pt 0\kern 1.5pt 20}=0.
\endaligned
\mytag{2.17}
$$
Substituting \mythetag{2.15} into \mythetag{2.17}, we derive the equation
$$
\hskip -2em
\aligned
{C_i}^{10}&+(6\,a_{22}^4-2\,a_{12}^4-a_{12}^2\,a_{22}^2)\,{C_i}^8
+(10\,a_{12}^2\,a_{22}^6+a_{12}^8\,+\\
&+\,a_{22}^8+4\,a_{12}^4\,a_{22}^4-14\,a_{12}^6\,a_{22}^2)\,{C_i}^6
+(14\,a_{12}^4\,a_{22}^8\,-\\
&-\,4\,a_{12}^6\,a_{22}^6-a_{12}^2\,a_{22}^{10}-10\,a_{12}^8\,a_{22}^4
-a_{12}^{10}\,a_{22}^2)\,{C_i}^4\,+\\
&+\,(2\,a_{12}^6\,a_{22}^{10}-6\,a_{12}^{10}\,a_{22}^6
+a_{12}^8\,a_{22}^8)\,{C_i}^2-a_{12}^{10}\,a_{22}^{10}=0.
\endaligned
\mytag{2.18}
$$
It is remarkable that the equation \mythetag{2.18} can be produced from the
equation \mythetag{1.1} by substituting $t=C_i$, $p=a_{12}$ and $q=a_{22}$.
Therefore the only known case where the equation \mythetag{2.18} has a rational
solution for $C_i$ is the case $a_{12}=a_{22}$. Relying on the irreducibility
of the fraction in \mythetag{2.3}, we set 
$$
\xalignat 2
&\hskip -2em
a_{12}=1,
&&a_{22}=1.
\mytag{2.19}
\endxalignat
$$
Then, in order to satisfy the equality \mythetag{2.4}, we choose
$$
\xalignat 2
&\hskip -2em
a_{11}=1,
&&a_{21}=0.
\mytag{2.20}
\endxalignat
$$
The choice \mythetag{2.20} is not unique. But it is the most simple.\par
     Applying \mythetag{2.19} and \mythetag{2.20} to \mythetag{2.7}, we derive
the following formulas: 
$$
\pagebreak
\xalignat 3
&\hskip -2em
\tilde p=p-q,
&&\tilde q=q,
&&p=\tilde p+\tilde q.
\quad
\mytag{2.21}
\endxalignat 
$$
Due to \mythetag{2.10} the formulas \mythetag{2.21} mean that we choose 
the bisectorial direction on the $p\,q$\,-\,coordinate plane for
asymptotic expansions.\par
\head
3. Bisectorial expansions. 
\endhead
     Let's substitute \mythetag{2.19} and \mythetag{2.20} into the equation
\mythetag{2.11}. As a result we can write the polynomial 
$Q_{\tilde p\kern 0.6pt \tilde q}(t)$ from \mythetag{2.11} in an explicit 
form:  
$$
\hskip -2em
\gathered
\hskip -2em
Q_{\tilde p\kern 0.6pt \tilde q}(t)=t^{10}
+(3\,\tilde q^{\kern 1pt 4}-10\,\tilde q^{\kern 1pt 3}\,\tilde p
-13\,\tilde q^{\kern 1pt 2}\,\tilde p^{\kern 1pt 2}
-8\,\tilde p^{\kern 1pt 3}\,\tilde q-2\,\tilde p^{\kern 1pt 4})\,t^8
+(\tilde p^{\kern 1pt 8}\,-\\
-\,40\,\tilde q^{\kern 1pt 7}\,\tilde p
-136\,\tilde q^{\kern 1pt 4}\,\tilde p^{\kern 1pt 4}+2\,\tilde q^{\kern 1pt 8}
+14\,\tilde p^{\kern 1pt 6}\,\tilde q^{\kern 1pt 2}
-148\,\tilde q^{\kern 1pt 6}\,\tilde p^{\kern 1pt 2}
-208\,\tilde q^{\kern 1pt 5}\,\tilde p^{\kern 1pt 3}\,+\\
+\,8\,\tilde p^{\kern 1pt 7}\,\tilde q
-28\,\tilde p^{\kern 1pt 5}\,\tilde q^{\kern 1pt 3})\,t^6
-(10\,\tilde p^{\kern 1pt 9}\,\tilde q^{\kern 1pt 3}
+302\,\tilde q^{\kern 1pt 10}\,\tilde p^{\kern 1pt 2}
+836\,\tilde q^{\kern 1pt 7}\,\tilde p^{\kern 1pt 5}\,+\\
+\,60\,\tilde q^{\kern 1pt 11}\,\tilde p
+200\,\tilde p^{\kern 1pt 7}\,\tilde q^{\kern 1pt 5}
+704\,\tilde q^{\kern 1pt 9}\,\tilde p^{\kern 1pt 3}
+\tilde p^{\kern 1pt 10}\,\tilde q^{\kern 1pt 2}
+494\,\tilde q^{\kern 1pt 6}\,\tilde p^{\kern 1pt 6}
+2\,\tilde q^{\kern 1pt 12}\,+\\
+\,956\,\tilde q^{\kern 1pt 8}\,\tilde p^{\kern 1pt 4}
+55\,\tilde p^{\kern 1pt 8}\,\tilde q^{\kern 1pt 4})\,t^4
-(6\,\tilde p^{\kern 1pt 10}\,\tilde q^{\kern 1pt 6}
+1444\,\tilde q^{\kern 1pt 11}\,\tilde p^{\kern 1pt 5}
+1230\,\tilde q^{\kern 1pt 10}\,\tilde p^{\kern 1pt 6}\,+\\
+\,712\,\tilde q^{\kern 1pt 9}\,\tilde p^{\kern 1pt 7}
+40\,\tilde q^{\kern 1pt 15}\,\tilde p
+1160\,\tilde q^{\kern 1pt 12}\,\tilde p^{\kern 1pt 4}
+624\,\tilde q^{\kern 1pt 13}\,\tilde p^{\kern 1pt 3}
+269\,\tilde q^{\kern 1pt 8}\,\tilde p^{\kern 1pt 8}\,+\\
+\,3\,\tilde q^{\kern 1pt 16}
+60\,\tilde p^{\kern 1pt 9}\,\tilde q^{\kern 1pt 7}
+212\,\tilde q^{\kern 1pt 14}\,\tilde p^{\kern 1pt 2})\,t^2
-210\,\tilde q^{\kern 1pt 16}\,\tilde p^{\kern 1pt 4}
-\tilde q^{\kern 1pt 10}\,\tilde p^{\kern 1pt 10}
-\tilde q^{\kern 1pt 20}\,-\\
-\,210\,\tilde q^{\kern 1pt 14}\,\tilde p^{\kern 1pt 6}
-120\,\tilde q^{\kern 1pt 13}\,\tilde p^{\kern 1pt 7}
-45\,\tilde q^{\kern 1pt 12}\,\tilde p^{\kern 1pt 8}
-10\,\tilde q^{\kern 1pt 11}\,\tilde p^{\kern 1pt 9}
-10\,\tilde q^{\kern 1pt 19}\,\tilde p\,\,-\\
-\,120\,\tilde q^{\kern 1pt 17}\,\tilde p^{\kern 1pt 3}
-252\,\tilde q^{\kern 1pt 15}\,\tilde p^{\kern 1pt 5}
-45\,\tilde q^{\kern 1pt 18}\,\tilde p^{\kern 1pt 2}.
\endgathered
\mytag{3.1}
$$
Now let's substitute \mythetag{2.19} into the equation \mythetag{2.18}. 
As a result we can factor it:
$$
\hskip -2em
({C_i}-1)\,({C_i}+1)\,({C_i}^2+1)^4=0.
\mytag{3.2}
$$ 
The equation \mythetag{3.2} has two simple real roots
$$
\xalignat 2
&C_i=1, &&C_i=-1
\endxalignat
$$
and two purely imaginary roots 
$$
\xalignat 2
&C_i=\goth i, &&C_i=-\goth i
\endxalignat
$$
of multiplicity four. Here $\goth i=\sqrt{-1}$. The condition 
\mythetag{2.12} excludes two roots $C_i=-1$ and $C_i=-\goth i$. 
Therefore from \mythetag{2.14} we derive the expansion 
$$
\hskip -2em
t_i(\tilde p,\tilde q)=\tilde q^{\kern 1pt 2}
\biggl(1+\sum^\infty_{s=1}\beta_{is}\,\tilde q^{-s}\biggr)
\text{\ \ as \ }\tilde q\to+\infty
\mytag{3.3}
$$
for real roots of the polynomial \mythetag{3.1} and the expansion 
$$
\hskip -2em
t_i(\tilde p,\tilde q)=\goth i\,\tilde q^{\kern 1pt 2}
\biggl(1+\sum^\infty_{s=1}\beta_{is}\,\tilde q^{-s}\biggr)
\text{\ \ as \ }\tilde q\to+\infty
\mytag{3.4}
$$
for complex roots of the polynomial \mythetag{3.1} in the equation
\mythetag{2.11}.\par
\head
4. Asymptotic estimates for real roots.
\endhead
     According to the formulas \mythetag{3.3} and \mythetag{3.4},
both real and complex roots of the polynomial \mythetag{3.1} are
growing as $\tilde q\to+\infty$. \pagebreak Therefore we need to specify the
expansions \mythetag{3.3} and \mythetag{3.4} up to non-growing
terms. For a real root we have 
$$
\hskip -2em
t_1(\tilde p,\tilde q)=\tilde q^{\kern 1pt 2}+
5\,\tilde p\,\tilde q+10\,\tilde p^{\kern 1pt 2}+R_1{(\tilde p,\tilde q)}
\text{\ \ as \ }\tilde q\to+\infty.
\mytag{4.1}
$$
The formula \mythetag{4.1} is in agreement with the formula \mythetag{3.3}. 
It means that $\beta_{11}=5\,\tilde p$ and $\beta_{12}=10\,\tilde p^2$. Like
in \mycite{1} and \mycite{2}, we have to obtain an estimate of the form
$$
\hskip -2em
|R_1{(\tilde p,\tilde q)}|<\frac{C(\tilde p)}{\tilde q}
\mytag{4.2}
$$
for the remainder term $R_1{(\tilde p,\tilde q)}$ in \mythetag{4.1}. For
this purpose we substitute 
$$
\hskip -2em
t=\tilde q^{\kern 1pt 2}+
5\,\tilde p\,\tilde q+10\,\tilde p^{\kern 1pt 2}
+\frac{c}{\tilde q}
\mytag{4.3}
$$
into the polynomial \mythetag{3.1}. Then we replace $\tilde q$ with the new
variable $z$:
$$
\hskip -2em
z=\frac{1}{\tilde q}.
\mytag{4.4}
$$
As a result of two substitutions \mythetag{4.3} and \mythetag{4.4} and upon 
removing denominators the equation \mythetag{2.11} with the polynomial
\mythetag{3.1} turns to a polynomial equation in the new variables $c$ and
$z$. This equation can be written as
$$
\hskip -2em
1216\,\tilde p^{\kern 1pt 3}+f(\tilde p,c,z)=-32\,c. 
\mytag{4.5}
$$
Here $f(\tilde p,c,z)$ is a polynomial of three variables given by an explicit 
formula. The formula for $f(\tilde p,c,z)$ is rather huge. Therefore it is 
placed to the ancillary file \darkred{{\tt strategy\kern -0.5pt\_\kern 1.5pt 
formulas\_\kern 0.5pt 03.txt}} in a machine-readable form.\par
     Let $\tilde q\geqslant 97\,|\tilde p\kern 0.5pt|$ and let $c$ belong to
one of the following two intervals:
$$
\hskip -2em
\aligned
-74\,|\tilde p\kern 0.5pt|^3<\ &c< 0\text{\ \ if \ }\tilde p>0,\\
0<\ &c< 74\,|\tilde p\kern 0.52pt|^3\text{\ \ if \ }\tilde p<0.
\endaligned
\mytag{4.6}
$$
From $\tilde q\geqslant 97\,|\tilde p\kern 0.5pt|$ and from \mythetag{4.4}
we derive the estimate $|z|\leqslant 1/97\,|\tilde p\kern 0.5pt|^{-1}$. Using
this estimate and using the inequalities \mythetag{4.6}, by means of direct
calculations one can derive the following estimate for the modulus of the 
function $f(\tilde p,c,z)$:	
$$
\hskip -2em
|f(\tilde p,c,z)|<1142\,|\tilde p\kern 0.5pt|^{3}.
\mytag{4.7}
$$
For fixed $\tilde p$ and $z$ the estimate \mythetag{4.7} means that the 
left hand side of the equation \mythetag{4.5} is a continuous function of 
$c$ whose values obey the inequalities 
$$
\hskip -2em
\aligned
74\,|\tilde p\kern 0.5pt|^3\leqslant\ &1216\,\tilde p^{\kern 1pt 3}
+f(\tilde p,c,z)\leqslant 2358\,\,|\tilde p\kern 0.5pt|^3\text{\ \ if \ }
\tilde p>0,\\
-2358\,|\tilde p\kern 0.5pt|^3\leqslant\ &1216\,\tilde p^{\kern 1pt 3}
+f(\tilde p,c,z)\leqslant -74\,|\tilde p\kern 0.5pt|^3\text{\ \ if \ }
\tilde p<0\\
\endaligned
\mytag{4.8}
$$
while $c$ runs over the corresponding interval \mythetag{4.6}. The right 
hand side of the equation \mythetag{4.5} is also a continuous function of 
the variable $c$. Moreover, it is monotonic. Multiplying the inequalities 
\mythetag{4.6} by $-32$, we find that the values of the right hand side 
of the equation \mythetag{4.5} fill one of the two intervals 
$$
\hskip -2em
\aligned
0<\ &-32\,c<2368\,|\tilde p\kern 0.5pt|^3 \text{\ \ if \ }\tilde p>0,\\
-2368\,|\tilde p\kern 0.5pt|^3<\ &-32\,c< 0\text{\ \ if \ }\tilde p<0
\endaligned
\mytag{4.9}
$$
while $c$ runs over the corresponding interval \mythetag{4.6}. Comparing 
\mythetag{4.8} with \mythetag{4.9}, we see that there is at least one root 
of the polynomial equation \mythetag{4.5} somewhere in one of the two
intervals \mythetag{4.6}.\par
     The case $\tilde p=0$ is exceptional. In this case $f(\tilde p,c,z)=0$.
Hence $c=0$ is a root of the equation \mythetag{4.5} for this case.\par
     The variable $c$ is related to the variable $t$ by means of the formula 
\mythetag{4.4}. Therefore the inequalities \mythetag{4.6} for $c$ imply the 
following inequalities for $t$:
$$
\hskip -2em
\gathered
\tilde q^{\kern 1pt 2}
+5\,\tilde p\,\tilde q+10\,\tilde p^{\kern 1pt 2}
-\frac{74\,|\tilde p\kern 0.5pt|^3}{\tilde q}<t
<\tilde q^{\kern 1pt 2}
+5\,\tilde p\,\tilde q+10\,\tilde p^{\kern 1pt 2}
\text{\ \ if \ }\tilde p>0,\\
\tilde q^{\kern 1pt 2}+
5\,\tilde p\,\tilde q+10\,\tilde p^{\kern 1pt 2}
<t<\tilde q^{\kern 1pt 2}
+5\,\tilde p\,\tilde q+10\,\tilde p^{\kern 1pt 2}
+\frac{74\,|\tilde p\kern 0.5pt|^3}{\tilde q}
\text{\ \ if \ }\tilde p<0.
\endgathered
\mytag{4.10}
$$
The case $\tilde p=0$ is exceptional. In this case we get
$$
\hskip -2em
t=\tilde q^{\kern 1pt 2}\text{\ \ if \ }\tilde p=0.
\mytag{4.11}
$$
The result obtained is formulated as a theorem. 
\mytheorem{4.1} For each $\tilde q\geqslant 97\,|\tilde p\kern 0.5pt|$
there is at least one real root of the polynomial \mythetag{3.1} satisfying 
one of the conditions in \mythetag{4.10} and \mythetag{4.11}.
\endproclaim
     Theorem~\mythetheorem{4.1} proves the asymptotic expansion 
\mythetag{4.1} and provides the estimate of the form \mythetag{4.2} for 
the remainder term in it. 
\head
5. Asymptotics for complex roots.
\endhead
     For complex roots of the polynomial \mythetag{3.1} we have the formula 
specifying \mythetag{3.4}:
$$
\hskip -2em
t_i(\tilde p,\tilde q)=\goth i\,\tilde q^{\kern 1pt 2}+
u_i\,\tilde p\,\tilde q-\frac{u_i+\goth i\,u_i^2}{2}
\,\tilde p^{\kern 1pt 2}+R_i{(\tilde p,\tilde q)}
\text{\ \ as \ }\tilde q\to+\infty
\mytag{5.1}
$$
and $i=2,\,\ldots,\,5$. Here $\goth i=\sqrt{-1}$ and $u_i$ are roots of 
the following quartic equation:
$$
\hskip -2em
u^4+8\,u^2-12\,\goth i\,u -4=0.
\mytag{5.2}
$$
The equation \mythetag{5.2} is irreducible. All of its roots are irrational. 
Two of them are purely imaginary. Here are approximate values for these
two roots:
$$
\xalignat 2
&\hskip -2em
u_2\approx 0.4863801704\,\goth i,
&&u_3\approx -3.439109107\,\goth i.
\mytag{5.3}
\endxalignat
$$ 
The other two roots are complex. Their approximate values are 
$$
\hskip -2em
\aligned
&u_4\approx\hphantom{-}0.4600767354+1.476364468\,\goth i,\\
&u_5\approx -0.4600767354+1.476364468\,\goth i.
\endaligned
\mytag{5.4}
$$\par
     Like in the previous section, below we derive an estimate of the form
$$
\hskip -2em
|R_i{(\tilde p,\tilde q)}|<\frac{C_i(\tilde p)}{\tilde q}
\mytag{5.5}
$$
for the remainder term in the formula \mythetag{5.1}. For
this purpose we substitute 
$$
\hskip -2em
t=\goth i\,\tilde q^{\kern 1pt 2}+
u\,\tilde p\,\tilde q-\frac{u+\goth i\,u^2}{2}
\,\tilde p^{\kern 1pt 2}+\frac{c}{\tilde q}.
\mytag{5.6}
$$
into the polynomial \mythetag{3.1}. Then we replace $\tilde q$ with the new
variable $z$ using \mythetag{4.4}. As a result of two substitutions \mythetag{5.6} 
and \mythetag{4.4} and upon removing denominators the equation \mythetag{2.11} 
with the polynomial \mythetag{3.1} turns to a polynomial equation in the new 
variables $c$ and $z$. This equation can be written as
$$
\hskip -2em
\aligned
212992\,\,\goth i\,\,\tilde p^{\kern 1pt 6}
&-598016\,\tilde p^{\kern 1pt 6}\,u
-446464\,\,\goth i\,\,\tilde p^{\kern 1pt 6}\,u^2\,+\\
&+\,110592\,\tilde p^{\kern 1pt 6}\,u^3
+\varphi(u,\tilde p,c,z)=352256\,\tilde p^{\kern 1pt 3}\,c.
\endaligned
\mytag{5.7}
$$
Here $\varphi(u,\tilde p,c,z)$ is a polynomial of four variables given by an 
explicit formula. The formula for $\varphi(u,\tilde p,c,z)$ is rather huge. 
Therefore it is placed to the ancillary file \darkred{{\tt strategy\kern -0.5pt
\_\kern 1.5pt formulas\_\kern 0.5pt 03.txt}} in a machine-readable form.\par
     Let $\tilde q\geqslant 97\,|\tilde p\kern 0.5pt|$. Now $c$ is a complex 
variable. Assume that it runs over the open disk on the complex plane given the 
by inequality 
$$
\hskip -2em
\aligned
|c|<51\,|\tilde p\kern 0.5pt|^3.
\endaligned
\mytag{5.8}
$$
From $\tilde q\geqslant 97\,|\tilde p\kern 0.5pt|$ and from \mythetag{4.4}
we derive the estimate $|z|\leqslant 1/97\,|\tilde p\kern 0.5pt|^{-1}$. Using
this estimate and using \mythetag{5.8}, upon substituting the value $u=u_2$ 
from \mythetag{5.3} one can derive the following estimate for the modulus of 
the function $\varphi(u,\tilde p,c,z)$:	
$$
\hskip -2em
|\varphi(u_2,\tilde p,c,z)|\leqslant 1174818\,|\tilde p\kern 0.5pt|^{6}.
\mytag{5.9}
$$
For fixed $\tilde p$ and $z$ the estimate \mythetag{5.9} means that the 
left hand side of the equation \mythetag{5.7} is a function of 
$c$ whose values within the disc \mythetag{5.8} obey the estimate
$$
\hskip -2em
\aligned
|212992\,\,\goth i\,\,\tilde p^{\kern 1pt 6}
&-598016\,\tilde p^{\kern 1pt 6}\,u_2
-446464\,\,\goth i\,\,\tilde p^{\kern 1pt 6}\,u_2^2\,+\\
&+\,110592\,\tilde p^{\kern 1pt 6}\,u_2^3
+\varphi(u_2,\tilde p,c,z)|\leqslant 1189840\,|\tilde p\kern 0.5pt|^{6}.
\endaligned
\mytag{5.10}
$$
Note that it is a holomorphic function which is continuous up to the boundary
of the disc \mythetag{5.8}. Therefore the estimate \mythetag{5.10} holds
on the boundary of the disc.\par
    The right hand side of the equation \mythetag{5.7} is also a holomorphic 
function of $c$. It has exactly one zero at the origin within the disc 
\mythetag{5.8} and its modulus is constant on the boundary of this disc.
Indeed, we have 
$$
\hskip -2em
|352256\,\tilde p^{\kern 1pt 3}\,c|=17965056\,|\tilde p\kern 0.5pt|^{6}
\text{\ \ if \ }|c|=74\,|\tilde p\kern 0.5pt|^3.
\mytag{5.11}
$$
Comparing the numbers $1189839<17965056$ from \mythetag{5.10} and \mythetag{5.11}
and applying the Rouch\'e theorem (see \mycite{68} or \mycite{69}), we conclude
that the equation \mythetag{5.7} with $u=u_2$ has exactly one solution within the 
disc \mythetag{5.8}. \par
     The same conclusion is valid for the other three roots $u=u_3$, $u=u_4$, 
and $u=u_5$ of the equation \mythetag{5.2} from \mythetag{5.3} and \mythetag{5.4}. 
However, the estimate \mythetag{5.10} for them is replaced by the following 
three estimates:
$$
\align
&\hskip -2em
\aligned
|212992\,\,\goth i\,\,\tilde p^{\kern 1pt 6}
&-598016\,\tilde p^{\kern 1pt 6}\,u_3
-446464\,\,\goth i\,\,\tilde p^{\kern 1pt 6}\,u_3^2\,+\\
&+\,110592\,\tilde p^{\kern 1pt 6}\,u_3^3
+\varphi(u_3,\tilde p,c,z)|\leqslant 
16504669\,|\tilde p\kern 0.5pt|^{6},
\endaligned
\mytag{5.12}\\
\vspace{1ex}
&\hskip -2em
\aligned
|212992\,\,\goth i\,\,\tilde p^{\kern 1pt 6}
&-598016\,\tilde p^{\kern 1pt 6}\,u_4
-446464\,\,\goth i\,\,\tilde p^{\kern 1pt 6}\,u_4^2\,+\\
&+\,110592\,\tilde p^{\kern 1pt 6}\,u_4^3
+\varphi(u_4,\tilde p,c,z)|\leqslant 
2513770\,|\tilde p\kern 0.5pt|^{6},
\endaligned
\mytag{5.13}\\
\vspace{1ex}
&\hskip -2em
\aligned
|212992\,\,\goth i\,\,\tilde p^{\kern 1pt 6}
&-598016\,\tilde p^{\kern 1pt 6}\,u_5
-446464\,\,\goth i\,\,\tilde p^{\kern 1pt 6}\,u_5^2\,+\\
&+\,110592\,\tilde p^{\kern 1pt 6}\,u_5^3
+\varphi(u_5,\tilde p,c,z)|\leqslant 
2513770\,|\tilde p\kern 0.5pt|^{6}.
\endaligned
\mytag{5.14}
\endalign
$$
The numbers $16504669$ and $2513770$ in the right hand sides of
\mythetag{5.12}, \mythetag{5.13}, and \mythetag{5.14} are smaller 
than the number $17965056$ in \mythetag{5.11}, which is the reason
for applying the Rouch\'e theorem.\par  
     The variable $c$ is related to the variable $t$ by means of the 
formula \mythetag{5.6}. Therefore the inequality \mythetag{5.8} for 
$c$ implies the following inequalities for $t$:
$$
\gather
\hskip -2em
\Bigl|\kern 1pt t-\goth i\,\tilde q^{\kern 1pt 2}
-u_2\,\tilde p\,\tilde q+\frac{u_2+\goth i\,u_2^2}{2}
\,\tilde p^{\kern 1pt 2}\Bigr|<\frac{51\,|\tilde p\kern 0.5pt|^{3}}
{\tilde q}\text{\ \ if \ }\tilde p\neq 0,
\mytag{5.15}\\
\vspace{1ex}
\hskip -2em
\Bigl|\kern 1pt t-\goth i\,\tilde q^{\kern 1pt 2}
-u_3\,\tilde p\,\tilde q+\frac{u_3+\goth i\,u_3^2}{2}
\,\tilde p^{\kern 1pt 2}\Bigr|<\frac{51\,|\tilde p\kern 0.5pt|^{3}}
{\tilde q}\text{\ \ if \ }\tilde p\neq 0,
\mytag{5.16}\\
\vspace{1ex}
\hskip -2em
\Bigl|\kern 1pt t-\goth i\,\tilde q^{\kern 1pt 2}
-u_4\,\tilde p\,\tilde q+\frac{u_4+\goth i\,u_4^2}{2}
\,\tilde p^{\kern 1pt 2}\Bigr|<\frac{51\,|\tilde p\kern 0.5pt|^{3}}
{\tilde q}\text{\ \ if \ }\tilde p\neq 0,
\mytag{5.17}\\
\vspace{1ex}
\hskip -2em
\Bigl|\kern 1pt t-\goth i\,\tilde q^{\kern 1pt 2}
-u_5\,\tilde p\,\tilde q+\frac{u_5+\goth i\,u_5^2}{2}
\,\tilde p^{\kern 1pt 2}\Bigr|<\frac{51\,|\tilde p\kern 0.5pt|^{3}}
{\tilde q}\text{\ \ if \ }\tilde p\neq 0.
\mytag{5.18}
\endgather
$$
The inequalities \mythetag{5.15}, \mythetag{5.16}, \mythetag{5.17}, and 
\mythetag{5.18} define four disks which play the same role as the intervals
\mythetag{4.10} in the previous section.\par
     The case $\tilde p=0$ is exceptional. In this case the disks 
\mythetag{5.15}, \mythetag{5.16}, \mythetag{5.17}, and \mythetag{5.18}
collapse to the point $t=\goth i\,\tilde q^{\kern 1pt 2}$ thus producing
a multiple root:
$$
\hskip -2em
t=\goth i\,\tilde q^{\kern 1pt 2}.
\mytag{5.19}
$$
Now we can formulate the result of this section as a theorem.
\mytheorem{5.1} For each $\tilde q\geqslant 97\,|\tilde p\kern 0.5pt|$
there is exactly one root of the polynomial \mythetag{3.1} in each of 
the four disks \mythetag{5.15}, \mythetag{5.16}, \mythetag{5.17}, 
and \mythetag{5.18} or there is one multiple root given by the
formula \mythetag{5.19}.
\endproclaim
     Theorem~\mythetheorem{5.1} proves the asymptotic expansion 
\mythetag{5.1} and provides the estimate of the form \mythetag{5.5} for 
the remainder terms in it. 
\head
6. Non-intersection of asymptotic sites.
\endhead
    In the previous two sections we have found five sites where 
roots of the polynomial \mythetag{3.1} are located. \pagebreak 
They are the intervals \mythetag{4.10} and the disks \mythetag{5.15}, 
\mythetag{5.16}, \mythetag{5.17}, \mythetag{5.18} in the non-degenerate 
case $\tilde p\neq 0$ and the points \mythetag{4.10} and \mythetag{5.19} 
in the degenerate case $\tilde p=0$. The points \mythetag{4.10} and 
\mythetag{5.19} do not coincide since $\tilde q=q>0$ (see \mythetag{2.21} 
and Theorem~\mythetheorem{1.1}). In the case $\tilde p\neq 0$ we have
the following lemma. 
\mylemma{6.1} For $\tilde q\geqslant 97\,|\tilde p\kern 0.5pt|\neq 0$
the asymptotic sites \mythetag{4.10}, \mythetag{5.15}, \mythetag{5.16}, 
\mythetag{5.17}, and \mythetag{5.18} do not intersect with each other. 
\endproclaim
\demo{Proof} In order to prove Lemma~\mythelemma{6.1} for discs it is 
sufficient to calculate the distances between their centers ans compare 
them with the double of their radius, which is the same for all of them. 
The intervals \mythetag{4.10} are centered about the point 
$t=\tilde q^{\kern 1pt 2}+5\,\tilde p\,\tilde q
+10\,\tilde p^{\kern 1pt 2}$. They can be enclosed in one disc with 
the radius $74\,|\tilde p\kern 0.5pt|^3/\tilde q$ for rough estimates. 
By means of direct calculations one can derive a lower bound
for the distances from the center of this disc to the centers of the 
other four discs: 
$$
\hskip -2em
d_{1\kern 0.4pt i}\geqslant 1.4\,(\tilde q-2.5\,|\tilde p\kern 0.5pt|)^2
-21.75\,|\tilde p\kern 0.5pt|^2,\quad i=2,\,\ldots,\,5.
\mytag{6.1}
$$
Applying $\tilde q\geqslant 97\,|\tilde p\kern 0.5pt|$ to 
\mythetag{6.1}, we obtain the inequality
$$
\hskip -2em
d_{1\kern 0.4pt i}\geqslant 12480\,|\tilde p\kern 0.5pt|^2.
\mytag{6.2}
$$
For the sum of the radii of two discs from $\tilde q\geqslant 97\,|\tilde p
\kern 0.5pt|$ we derive
$$
r_1+r_i=\frac{74\,|\tilde p\kern 0.5pt|^3}{\tilde q}
+\frac{51\,|\tilde p\kern 0.5pt|^{3}}{\tilde q}\leqslant
\frac{125}{97}\,|\tilde p\kern 0.5pt|^2<2\,\,|\tilde p\kern 0.5pt|^2.
\mytag{6.3}
$$
Comparing \mythetag{6.2} and \mythetag{6.3}, we see that 
$d_{1\kern 0.4pt i}>r_1+r_i$, i\.\,e\. the intervals \mythetag{4.10}
do not intersect with the discs \mythetag{5.15}, \mythetag{5.16}, 
\mythetag{5.17}, and \mythetag{5.18}.\par
     Similarly, for the the mutual distances between centers of the
discs \mythetag{5.15}, \mythetag{5.16}, \mythetag{5.17}, and \mythetag{5.18}
one can obtain the following lower bound:
$$
\hskip -2em
d_{ij}\geqslant 0.98\,|\tilde p\kern 0.5pt|\,\tilde q
-8\,|\tilde p\kern 0.5pt|^2,\quad i=2,\,\ldots,\,5\text{\ \ and \ }
i\neq j.
\mytag{6.4}
$$
Applying $\tilde q\geqslant 97\,|\tilde p\kern 0.5pt|$ to 
\mythetag{6.4}, we obtain the inequality
$$
\hskip -2em
d_{1\kern 0.4pt i}\geqslant 87\,|\tilde p\kern 0.5pt|^2.
\mytag{6.5}
$$
For the double radius of the discs \mythetag{5.15}, \mythetag{5.16}, 
\mythetag{5.17}, and \mythetag{5.18} from the inequality $\tilde q
\geqslant 97\,|\tilde p\kern 0.5pt|$ we derive the following upper
bound: 
$$
\hskip -2em
r_i+r_j=2\,r_i=\frac{102\,|\tilde p\kern 0.5pt|^{3}}{\tilde q}
\leqslant \frac{102}{97}\,\,|\tilde p\kern 0.5pt|^2
<2\,\,|\tilde p\kern 0.5pt|^2.
\mytag{6.6}
$$
Comparing \mythetag{6.5} and \mythetag{6.6}, we see that 
$d_{ij}>r_i+r_j$, i\.\,e\. the discs \mythetag{5.15}, \mythetag{5.16}, 
\mythetag{5.17}, and \mythetag{5.18} do not intersect with each 
other. Lemma~\mythelemma{6.1} is proved.
\qed\enddemo
\mylemma{6.2} For $\tilde q\geqslant 97\,|\tilde p\kern 0.5pt|\neq 0$
the asymptotic discs \mythetag{5.15}, \mythetag{5.16}, \mythetag{5.17}, 
and \mythetag{5.18} are enclosed in upper half-plane of the complex
plane and do not intersect with the real axis. 
\endproclaim
\demo{Proof} The proof is based on the following lower bound for the 
distances from the centers of the discs \mythetag{5.15}, \mythetag{5.16}, 
\mythetag{5.17}, \mythetag{5.18} to the real axis:
$$
\hskip -2em
d_i\geqslant (\tilde q-2\,|\tilde p\kern 0.5pt|)^2
-12\,|\tilde p\kern 0.5pt|^2,\quad i=2,\,\ldots,\,5.
\mytag{6.7}
$$
Applying $\tilde q\geqslant 97\,|\tilde p\kern 0.5pt|$ to 
\mythetag{6.7}, we obtain the inequality
$$
\hskip -2em
d_i\geqslant 9013\,|\tilde p\kern 0.5pt|^2.
\mytag{6.8}
$$
On the other hand, for the radius of the discs \mythetag{5.15}, 
\mythetag{5.16}, \mythetag{5.17}, and \mythetag{5.18} from the inequality 
$\tilde q \geqslant 97\,|\tilde p\kern 0.5pt|$ we derive the following 
upper bound: 
$$
\hskip -2em
r_i=\frac{51\,|\tilde p\kern 0.5pt|^{3}}{\tilde q}
\leqslant \frac{51}{97}\,\,|\tilde p\kern 0.5pt|^2
<|\tilde p\kern 0.5pt|^2.
\mytag{6.9}
$$
Comparing \mythetag{6.8} and \mythetag{6.9}, we see that 
Lemma~\mythelemma{6.2} is proved.
\qed\enddemo
     Lemmas~\mythelemma{6.1} and \mythelemma{6.2} are summed up in the 
following theorem. 
\mytheorem{6.1} For $\tilde q \geqslant 97\,|\tilde p\kern 0.5pt|\neq 0$
five roots $t_1$, $t_2$, $t_3$, $t_4$, $t_5$ of the polynomial 
\mythetag{3.1} obeying the condition \mythetag{2.12} are simple. They are 
located within five disjoint sites \mythetag{4.10}, \mythetag{5.15}, 
\mythetag{5.16}, \mythetag{5.17}, and \mythetag{5.18}, one per each 
site. 
\endproclaim
     Due to \mythetag{2.13} Theorem~\mythetheorem{6.1} locates all of the ten
roots of the polynomial \mythetag{3.1}. Theorem~\mythetheorem{6.1} does not
cover the degenerate case $\tilde p=0$. However, due to \mythetag{2.21} the 
equality $\tilde p=0$ implies $p=q$. Therefore the degenerate case $\tilde p=0$
does not produce perfect cuboids (see Theorem~\mythetheorem{1.1}).
\head
7. Integer points of asymptotic sites.
\endhead
      According to Theorem~\mythetheorem{6.1}, for $\tilde q \geqslant 97
\,|\tilde p\kern 0.5pt|\neq 0$ the equation \mythetag{2.11} with the polynomial
\mythetag{3.1} has exactly one real positive root $t_1$ belonging to one of the
two asymptotic intervals \mythetag{4.10}. The following theorem is immediate 
from \mythetag{4.10}.
\mytheorem{7.1} If\/ $\tilde q \geqslant 97\,|\tilde p\kern 0.5pt|\neq 0$ and
if\/ $\tilde q>74\,|\tilde p\kern 0.5pt|^3$, then the asymptotic intervals 
\mythetag{4.10} have no integer points. 
\endproclaim
\head
8. Application to the cuboid problem.
\endhead
     The equation \mythetag{1.1} is related to the perfect cuboid
problem through Theorem~\mythetheorem{1.1}. The equation \mythetag{2.11} 
differs from the equation \mythetag{1.1} by the change of variables 
\mythetag{2.21}. Let's consider the case $\tilde p<0$ in \mythetag{4.10}.
In this case from \mythetag{4.10} we take
$$
\hskip -2em
t<\tilde q^{\kern 1pt 2}
+5\,\tilde p\,\tilde q+10\,\tilde p^{\kern 1pt 2}
+\frac{74\,|\tilde p\kern 0.5pt|^3}{\tilde q}.
\mytag{8.1}
$$
Theorem~\mythetheorem{1.1} provides four additional inequalities 
\mythetag{1.3}. Since $q=\tilde q$ in \mythetag{2.21}, the third of 
them is written as $t>\tilde q^{\kern 1pt 2}$. Combining it with
\mythetag{8.1}, we get
$$
\hskip -2em
\tilde q^{\kern 1pt 2}
+5\,\tilde p\,\tilde q+10\,\tilde p^{\kern 1pt 2}
+\frac{74\,|\tilde p\kern 0.5pt|^3}{\tilde q}>
\tilde q^{\kern 1pt 2}.
\mytag{8.2}
$$
Since $\tilde p<0$, the inequality \mythetag{8.2} turns to the following
one:
$$
\hskip -2em
-5\,|\tilde p\kern 0.5pt|\,\tilde q+10\,|\tilde p\kern 0.5pt|^2
+\frac{74\,|\tilde p\kern 0.5pt|^3}{\tilde q}>0.
\mytag{8.3}
$$
Let's apply the inequality $\tilde q \geqslant 97\,|\tilde p\kern 0.5pt|$ 
from Theorem~\mythetheorem{7.1}. Since $\tilde q=q>0$, it yields
$$
\xalignat 2
&\hskip -2em
-5\,|\tilde p\kern 0.5pt|\,\tilde q\leqslant 
-5\cdot 97\,|\tilde p\kern 0.5pt|^2,
&&\frac{74\,|\tilde p\kern 0.5pt|^3}{\tilde q}\leqslant 
\frac{74\,|\tilde p\kern 0.5pt|^2}{97}.
\mytag{8.4}
\endxalignat
$$
From \mythetag{8.4} one can easily derive the following inequality:
$$
\hskip -0.01em
-5\,|\tilde p\kern 0.5pt|\,\tilde q+10\,|\tilde p\kern 0.5pt|^2
+\frac{74\,|\tilde p\kern 0.5pt|^3}{\tilde q}\leqslant
\Bigl(-5\cdot 97+10+\frac{74}{97}\Bigr)\,|\tilde p\kern 0.5pt|^2=
-\frac{46001}{97}\,|\tilde p\kern 0.5pt|^2<0.
\mytag{8.5}
$$
The inequalities \mythetag{8.3} and \mythetag{8.5} contradict each other.
The contradiction obtained means that Theorem~\mythetheorem{7.1} can be
modified in the following way. 
\mytheorem{8.1} If\/ $\tilde q \geqslant 97\,|\tilde p\kern 0.5pt|$ and\/
$\tilde p<0$, then the corresponding asymptotic interval in \mythetag{4.10} 
has no integer points producing perfect cuboids. 
\endproclaim
     Unfortunately Theorem~\mythetheorem{8.1} cannot be extended to the case
$\tilde p>0$. In this case the inequality $\tilde q \geqslant 97\,|\tilde p
\kern 0.5pt|$ does not contradict the cuboid inequalities \mythetag{1.3}. 
However, Theorem~\mythetheorem{7.1} is still valid in the case $\tilde p>0$. 
\par
Theorem~\mythetheorem{7.1} provides two inequalities $\tilde q \geqslant 97
\,|\tilde p\kern 0.5pt|\neq 0$ and $\tilde q>74\,|\tilde p\kern 0.5pt|^3$. 
In the case $\tilde p>0$, passing to the original variables $p$ and $q$ by 
means of the formulas \mythetag{2.21}, these two inequalities are written in 
the following way:
$$
\xalignat 2
&\hskip -2em
p\leqslant\frac{98}{97}\,q=q+\frac{q}{97},
&& p<q+\root3\of{\frac{q}{74}},
\mytag{8.6}
\endxalignat 
$$
Due to \mythetag{2.21} the inequality $\tilde p>0$ itself means $p>q$.\par
     Similarly, Theorem~\mythetheorem{8.1} provides the inequality $\tilde q 
\geqslant 97\,|\tilde p\kern 0.5pt|$ in the case $\tilde p<0$, i\.\,e\. 
when $p<q$. Passing to the original variables $p$ and $q$ by means of 
\mythetag{2.21}, we can write this inequality in the following way:
$$
\hskip -2em
p\geqslant\frac{96}{97}\,q=q-\frac{q}{97},
\mytag{8.7}
$$
Since the bisector line $p=q$ does not produce perfect cuboids, the 
inequalities \mythetag{8.6} and \mythetag{8.7} can be united and then 
written as follows:
$$
\xalignat 2
&\hskip -2em
q-\frac{q}{97}\leqslant p,
&&p\leqslant q+\min\Bigl(\frac{q}{97},\root3\of{\frac{q}{74}}\Bigr).
\mytag{8.8}
\endxalignat 
$$
In this form the above inequalities \mythetag{8.8} are similar to the 
linear inequalities \mythetag{1.9} defining the linear region.\par
\head
9. Conclusions. 
\endhead
     Theorems~\mythetheorem{7.1} and \mythetheorem{8.1} along with the
inequalities \mythetag{8.8} constitute the main result of this paper. 
The inequalities \mythetag{8.8} define a subregion within the linear 
region \mythetag{1.9}. Theorems~\mythetheorem{7.1} and \mythetheorem{8.1}
prove that no cuboids are available within this subregion.\par
     The subregion defined by the inequalities \mythetag{8.8} 
is rather small. It looks like a narrow spiky strip surrounding the 
bisector line $p=q$ within the positive quadrant of the 
$p\,q$\,-\,coordinate plane. A substantial part of the linear region
\mythetag{1.9} still remains for numeric search of perfect cuboids. 
\par 

\enddocument
\end